\documentclass{article}
\usepackage[affil-it]{authblk}
\usepackage[utf8]{inputenc}
\usepackage{amsmath, amsfonts, amsthm, amssymb, latexsym}
\usepackage{geometry}
\usepackage{enumerate}
\usepackage{multirow}
\usepackage{rotating}
\usepackage{graphicx,color}
\usepackage{float}
\usepackage{mathrsfs}

\newtheorem{thm}{Theorem}[section]

\newtheorem{prop}[thm]{Proposition}

\newtheorem*{remark}{Remark}

\title{Determining The Magnetic Potential In The Fractional Magnetic 
Calder\'on Problem}
\author{Li Li}
\affil{Department of Mathematics, University of Washington,\\
Seattle, WA 98195, USA}
\date{}

\begin{document}

\maketitle

\noindent \textbf{ABSTRACT.}\, We determine both the magnetic potential and the electric potential from the exterior partial measurements of the Dirichlet-to-Neumann map in the fractional linear magnetic Calder\'on problem by using an integral identity. We also determine both the magnetic potential and the nonlinearity in the fractional semilinear magnetic Calder\'on problem by using a first order linearization.

\section{Introduction}

In this paper we continue the study of the fractional magnetic Calder\'on problem introduced in \cite{li2020calderon}, which is a generalization of the fractional Calder\'on problem introduced in \cite{ghosh2017calderon, ghosh2020calderon} as well as a nonlocal analogue of the classical magnetic Calder\'on problem studied in \cite{ferreira2007determining, krupchyk2014uniqueness, nakamura1995global, sun1993inverse}. 

We consider the operator 
$\mathcal{R}^s_A$, which is formally defined by
\begin{equation}\label{Refml}
\mathcal{R}^s_A u(x):= 2\lim_{\epsilon\to 0^+}
\int_{\mathbb{R}^n\setminus B_\epsilon(x)}(u(x)- R_A(x, y)u(y))K(x,y)\,dy
\end{equation}
where $K$ is a function associated with a heat kernel (see Subsection 2.2 in \cite{li2020calderon}) satisfying
$$K(x, y)= K(y, x),\qquad K(x, y)\sim \frac{1}{|x-y|^{n+2s}},$$ 
$A$ is a real vector-valued magnetic potential and
\begin{equation}\label{cos}
R_A(x, y):= \cos((x-y)\cdot A(\frac{x+y}{2})).
\end{equation}
For real-valued $u$, $\mathcal{R}^s_A$ is the real part of the fractional operator $\mathcal{L}^s_A$ introduced in \cite{li2020calderon} and clearly we have $\mathcal{R}^s_A= \mathcal{R}^s_{-A}$.

Under appropriate assumptions on $A$ and the electric potential $q$, 
the exterior Dirichlet problem
\begin{equation}\label{Rsalep}
\left\{
\begin{aligned}
(\mathcal{R}^s_A+ q)u&= 0\quad \text{in}\,\,\Omega\\
u&= g\quad \text{in}\,\,\Omega_e\\
\end{aligned}
\right.
\end{equation}
is well-posed so we can define the solution operator $P_{A, q}: g\to u_g$
and the Dirichlet-to-Neumann map $\Lambda_{A, q}$, 
which is formally given by
\begin{equation}\label{LDN}
\Lambda_{A, q}g:= \mathcal{R}^s_A (P_{A, q}g)|_{\Omega_e}.
\end{equation}

In \cite{li2020calderon}, we determined $q$ from the exterior partial measurements of $\Lambda_{A, q}$
for fixed $A$. Our goal here is to determine both $A$ and $q$ from the knowledge of the Dirichlet-to-Neumann map. The following theorem is our first main result in this paper.

\begin{thm}
Suppose $\mathrm{supp}\,A_j\subset\Omega\subset B_r(0)$ for some constant $r> 0$,
$c\leq q_j\in L^\infty(\Omega)$ for some constant $c> 0$, 
$A_j\in L^\infty(\mathbb{R}^n)$, $W_j$ are open sets s.t. $W_j\cap B_{3r}(0)= \emptyset$
($j= 1, 2$). Let
$$W^{(1, 2)}= \{\frac{x+ y}{2}: x\in W_1, y\in W_2\}.$$
Also assume $W^{(1, 2)}\setminus (\mathrm{supp}\,A_1\cup \mathrm{supp}\,A_2)\neq \emptyset$. If 
\begin{equation}\label{LIDN}
\Lambda_{A_1, q_1}g|_{W_2}= \Lambda_{A_2, q_2}g|_{W_2}
\end{equation}
for any $g\in C^\infty_c(W_1)$, then $A_1= \pm A_2$ and $q_1= q_2$.
\end{thm}

\begin{remark}
Recall that in the classical magnetic Calder\'on problem, it is impossible to completely determine $A$ from the knowledge of the Dirichlet-to-Neumann map
since the Dirichlet-to-Neumann maps associated with $(A_1, q)$ and $(A_2, q)$
coincide whenever $A_1$ and $A_2$ are gauge equivalent, i.e.
$$A_1- A_2= \nabla\phi$$
for some smooth $\phi$ in $\bar{\Omega}$ with $\phi|_{\partial\Omega}= 0$.
However in this fractional magnetic Calder\'on problem, we are able to completely determine $A$ (up to the sign) from the knowledge of $\Lambda_{A, q}$. See \cite{covi2020inverse}
for the study of a different fractional magnetic Calder\'on problem and see \cite{covi2020some} for
a recent study connecting local and nonlocal magnetic Calder\'on problems.
More results for the fractional linear Calder\'on problem can be found in \cite{ghosh2020uniqueness, ruland2020fractional}.

\end{remark}

In this paper we also continue the study of the fractional semilinear magnetic Calder\'on problem introduced in \cite{li2020semilinear}, which is a generalization of the fractional semilinear Calder\'on problem introduced in 
\cite{lai2020inverse} as well as a semilinear analogue of the inverse problem studied in \cite{li2020calderon}. See \cite{isakov1994global, krupchyk2020remark, lassas2019inverse}
for results for the classical semilinear Calder\'on problem.

As in \cite{li2020semilinear}, we focus on $K(x, y)= c_{n, s}/|x-y|^{n+2s}$ in 
the fractional semilinear Calder\'on problem. In this case, $\mathcal{R}^s_A$ is the real part of the fractional magnetic Laplacian $(-\Delta)^s_A$ studied in \cite{d2018ground, squassina2016bourgain} for real-valued $u$.

We consider the nonlinear exterior Dirichlet problem
\begin{equation}\label{Rsanonlinear}
\left\{
\begin{aligned}
\mathcal{R}^s_A u+ a(x, u)&= 0\quad \text{in}\,\,\Omega\\
u&= g\quad \text{in}\,\,\Omega_e\\
\end{aligned}
\right.
\end{equation}
where the nonlinearity $a(x, z): \Omega\times \mathbb{R}\to \mathbb{R}$ satisfies\\
(i)\, $z\to a(\cdot, z)$ is analytic with values in 
the H\"older space $C^s(\Omega)$;\\
(ii)\, $a(x, 0)= 0$ and $\partial_z a(x, 0)\geq c> 0$
for some constant $c> 0$\\
so we have the Taylor's expansion
\begin{equation}\label{TaylorHolder}
a(x, z)= \sum^\infty_{k=1}a_k(x)\frac{z^k}{k!},\quad 
a_k(x)= \partial^k_z a(x, 0)\in C^s(\Omega)
\end{equation}
where the series converges in $C^s(\Omega)$ topology. 

Under some boundedness conditions on $A$, (\ref{Rsanonlinear}) is well-posed for small $g$ so we can define the solution operator $Q_{A, a}: g\to u_g$
and the Dirichlet-to-Neumann map $\Lambda_{A, a}$, 
which is formally given by
\begin{equation}\label{NDN}
\Lambda_{A, a}g:= \mathcal{R}^s_A (Q_{A, a}g)|_{\Omega_e}.
\end{equation}

In \cite{li2020semilinear}, we determined $a$ from the exterior partial measurements of $\Lambda_{A, a}$
for fixed $A$. Our next goal here is to determine both $A$ and $a$ from the knowledge of the Dirichlet-to-Neumann map. The following theorem is our second main result in this paper.
\begin{thm}
Suppose $\mathrm{supp}\,A_j\subset\Omega \subset B_r(0)$ for some constant $r> 0$
and $||A_j||_{L^\infty(\mathbb{R}^n)}\leq \pi/(8\sqrt{n}r)$, $a^{(j)}$ satisfy (i) and (ii), $W_j$ are open sets s.t. $W_j\cap B_{3r}(0)= \emptyset$ ($j= 1, 2$). Let
$$W^{(1, 2)}= \{\frac{x+ y}{2}: x\in W_1, y\in W_2\}.$$
Also assume $W^{(1, 2)}\setminus (\mathrm{supp}\,A_1\cup \mathrm{supp}\,A_2)\neq \emptyset$. If 
\begin{equation}\label{NIDN}
\Lambda_{A_1, a^{(1)}}g|_{W_2}= \Lambda_{A_2, a^{(2)}}g|_{W_2},\qquad g\in C^\infty_c(W_1)
\end{equation}
whenever $||g||_{C^2(\mathbb{R}^n)}$ is sufficiently small, then 
$A_1= \pm A_2$ in $\Omega$ and $a^{(1)}= a^{(2)}$ in $\Omega\times \mathbb{R}$.
\end{thm}

The rest of this paper is organized in the following way. In Section 2, we summarize the background knowledge. We prove Theorem 1.1 by using an integral identity in Section 3. Based on Theorem 1.1, we prove Theorem 1.2 by using a first order linearization in the Sobolev space $H^s(\mathbb{R}^n)$ in Section 4.~\\

\noindent \textbf{Acknowledgement.} The author is partly supported by National Science Foundation. The author would like to thank Professor Gunther Uhlmann for suggesting the problem and for helpful discussions.

\section{Preliminaries}

Throughout this paper
\begin{itemize}
\item $n\geq 2$ denotes the space dimension and
$0< s< 1$ denotes the fractional power

\item $\Omega$ denotes a bounded domain with $C^{1,1}$ boundary and
$\Omega_e:= \mathbb{R}^n\setminus\bar{\Omega}$

\item $B_r(0)$ denotes the open ball centered at the origin with radius $r> 0$

\item $A: \mathbb{R}^n\to \mathbb{R}^n$ denotes a real vector-valued magnetic potential

\item $c, C, C', C_1,\cdots$ denote positive constants (which may depend on some parameters but always independent of small constants $\epsilon, \rho$)

\item $\int\cdots\int= \int_{\mathbb{R}^n}\cdots\int_{\mathbb{R}^n}$

\item $X^*$ denotes the continuous dual space of $X$ and write
$\langle f, u\rangle= f(u)$ for $u\in X,\,f\in X^*$

\item $||\cdot||_{C^2(\mathbb{R}^n)}$ is defined by
$$||f||_{C^2(\mathbb{R}^n)}= 
\sum_{|\alpha|\leq 2}||\partial^\alpha f||_{L^\infty(\mathbb{R}^n)}.$$
\end{itemize}

\subsection{Function Spaces}

Throughout this paper we refer all function spaces to real-valued function spaces.

For $t\in \mathbb{R}$, 
$H^t(\mathbb{R}^n)$ denotes the Sobolev space $W^{t, 2}(\mathbb{R}^n)$.

We have the natural identification
$$H^{-t}(\mathbb{R}^n)= H^t(\mathbb{R}^n)^*.$$
Let $U$ be an open set and $F$ be a closed set in $\mathbb{R}^n$, 
$$H^t(U):= \{u|_U: u\in H^t(\mathbb{R}^n)\},\qquad 
H^t_F(\mathbb{R}^n):= \{u\in H^t(\mathbb{R}^n): \mathrm{supp}\,u\subset F\},$$
$$\tilde{H}^t(U):= 
\mathrm{the\,\,closure\,\,of}\,\, C^\infty_c(U)\,\,\mathrm{in}\,\, H^t(\mathbb{R}^n).$$
Since $\Omega$ is a bounded domain with $C^{1,1}$ boundary implies $\Omega$ is Lipschitz bounded, 
$$\tilde{H}^t(\Omega)= H^t_{\bar{\Omega}}(\mathbb{R}^n).$$

For $0< s< 1$, $C^s(U)$ denotes the H\"older space $C^{0, s}(U)$.

\subsection{Old Results}

All the results presented in this subsection can be found in \cite{li2020calderon} and \cite{li2020semilinear}.

We have the bilinear form definition of $\mathcal{R}^s_A$, which is given by
\begin{equation}\label{blrsa}
\langle \mathcal{R}^s_Au, v\rangle= 2\iint(u(x)-R_A(x,y)u(y))v(x)K(x,y)\,dxdy.
\end{equation}

\begin{prop}
Suppose $0< s< 1$ and $A\in L^\infty(\mathbb{R}^n)$, then the operator
$\mathcal{R}^s_A: H^s(\mathbb{R}^n)\to H^{-s}(\mathbb{R}^n)$
is linear, bounded and
\begin{equation}\label{Rsasym}
\langle \mathcal{R}^s_Au, v\rangle= \langle \mathcal{R}^s_Av, u\rangle.
\end{equation}
\end{prop}

We have the well-posedness of the linear and semilinear exterior Dirichlet problems.
\begin{prop}
Suppose $A\in L^\infty(\mathbb{R}^n)$ and $c\leq q\in L^\infty(\Omega)$,
then the linear exterior problem (\ref{Rsalep}) has a unique (weak) solution $u_g\in H^s(\mathbb{R}^n)$ for each $g\in H^s(\mathbb{R}^n)$ and the solution operator 
$P_{A, q}$ is bounded on $H^s(\mathbb{R}^n)$.
\end{prop}

\begin{prop}
Suppose $\mathrm{supp}\,A\subset\Omega \subset B_r(0)$ for some $r> 0$
and $||A||_{L^\infty(\mathbb{R}^n)}\leq \pi/(8\sqrt{n}r)$, $W\cap B_{3r}(0)= \emptyset$ and $g\in C^\infty_c(W)$. There exists a small
constant $\rho >0$ s.t. if $||g||_{C^2(\mathbb{R}^n)}\leq \rho$, 
then the semilinear exterior problem
(\ref{Rsanonlinear})
has a unique solution $u\in H^s(\mathbb{R}^n)\cap C^s(\mathbb{R}^n)$ satisfying
$$(u- P_{A, a_1}g)\in
M:=\{v\in C^s(\mathbb{R}^n): v|_{\Omega_e}= 0, ||v||_{C^s(\mathbb{R}^n)}\leq \rho\}.$$
Denote the associated solution operator by $Q_{A, a}$. Moreover, we have
$$||Q_{A, a}g||_{C^s(\mathbb{R}^n)}\leq C||g||_{C^2(\mathbb{R}^n)}.$$
\end{prop}

The Dirichlet-to-Neumann map $\Lambda_{A, q}$ associated with (\ref{Rsalep}) can be defined by a bilinear form on $H^s(\Omega_e)\times H^s(\Omega_e)$. For $g, h\in C^\infty_c(\Omega_e)$, we have 
\begin{equation}\label{bilinearDN}
\langle \Lambda_{A, q}g, h\rangle= \langle \mathcal{R}^s_A(P_{A, q}g), h^*\rangle
+ \int_\Omega (P_{A, q}g) h^*
\end{equation}
where $h^*\in H^s(\mathbb{R}^n)$ satisfying
$h^*- h\in \tilde{H}^s(\Omega)$.
We know that (\ref{bilinearDN}) does not depend on the choice of $h^*$ and
\begin{equation}\label{DNsym}
\langle \Lambda_{A, q}g, h\rangle= \langle \Lambda_{A, q}h, g\rangle.
\end{equation}
This bilinear form definition coincides with the definition given by (\ref{LDN}).

The Dirichlet-to-Neumann map $\Lambda_{A, a}$ associated with (\ref{Rsanonlinear})
cannot be defined by a bilinear form due to the nonlinearity $a(\cdot, \cdot)$. Proposition 2.3 ensures that (\ref{NDN}) is well-defined at least for $g$ satisfying the condition assumed in the statement of the proposition.

We also have the following Runge approximation property.
\begin{prop}
Suppose $\mathrm{supp}\,A\subset\Omega \subset B_r(0)$ for some $r> 0$, $W$ is an open set s.t. $W\subset \Omega_e$ and $W\cap B_{3r}(0)= \emptyset$, then
$$S:= \{P_{A, q}f|_\Omega: f\in C^\infty_c(W)\}$$
is dense in $L^2(\Omega)$.
\end{prop}

We will prove Theorem 1.1 based on the following theorem in \cite{li2020calderon}.
\begin{thm}
Suppose $\mathrm{supp}\,A\subset\Omega\subset B_r(0)$ for some constant $r> 0$,
$c\leq q_j\in L^\infty(\Omega)$ for some constant $c> 0$, 
$A\in L^\infty(\mathbb{R}^n)$, 
$W_j$ are open sets s.t. $W_j\cap B_{3r}(0)= \emptyset$
($j= 1, 2$). If 
$$\Lambda_{A, q_1}g|_{W_2}= \Lambda_{A, q_2}g|_{W_2}$$
for any $g\in C^\infty_c(W_1)$, then $q_1= q_2$.
\end{thm}

The following first order linearization relates the semilinear problem to the linear one.
\begin{prop}
Suppose $\mathrm{supp}\,A\subset\Omega \subset B_r(0)$ for some $r> 0$
and $||A||_{L^\infty(\mathbb{R}^n)}\leq \pi/(8\sqrt{n}r)$, $W\cap B_{3r}(0)= \emptyset$ and $g\in C^\infty_c(W)$, then $$Q_{A, a}(\epsilon g)/\epsilon\to P_{A, a_1}g$$ in $H^s(\mathbb{R}^n)$ as $\epsilon\to 0$.
\end{prop}

We will use the first order linearization above to prove Theorem 1.2 based on the following theorem in \cite{li2020semilinear}.
\begin{thm}
Suppose $\mathrm{supp}\,A\subset\Omega \subset B_r(0)$ for some constant $r> 0$
and $||A||_{L^\infty(\Omega)}\leq \pi/(8\sqrt{n}r)$, $a^{(j)}$ satisfy (i) and (ii), $W_j$ are open sets s.t. $W_j\cap B_{3r}(0)= \emptyset$ ($j= 1, 2$). If 
$$\Lambda_{A, a^{(1)}}g|_{W_2}= \Lambda_{A, a^{(2)}}g|_{W_2},\qquad g\in C^\infty_c(W_1)$$
whenever $||g||_{C^2(\mathbb{R}^n)}$ is sufficiently small, then $a^{(1)}= a^{(2)}$ in $\Omega\times \mathbb{R}$.
\end{thm}

\begin{remark}
Compared with the statement of Theorem 2.5 (respectively, Theorem 2.7), the statement of Theorem 1.1 (respectively, Theorem 1.2) contains an additional assumption on the set $W^{(1, 2)}$. This assumption ensures that the double integral
in the integral identity (obtained in Section 3) is actually over the region $\Omega\times \Omega$, which enables us to apply the Runge approximation property (Proposition 2.4). See the proof in Section 3 for details.
\end{remark}

\section{The Proof of Theorem 1.1}
We first build an integral identity, which will be useful in
the proof of Theorem 1.1.

For $g_j\in C^\infty_c(\Omega_e)$ and $u_j= P_{A_j, q_j}g_j$ ($j= 1,2$) solving
\begin{equation*}
\left\{
\begin{aligned}
(\mathcal{R}^s_{A_j}+ q_j)u&= 0\quad \text{in}\,\,\Omega\\
u&= g_j\quad \text{in}\,\,\Omega_e,\\
\end{aligned}
\right.
\end{equation*}
by (\ref{blrsa}), (\ref{Rsasym}), (\ref{bilinearDN}) and (\ref{DNsym}) we have
$$\langle (\Lambda_{A_1, q_1} -\Lambda_{A_2, q_2})g_1, g_2\rangle= 
\langle \Lambda_{A_1, q_1}g_1, g_2\rangle- \langle \Lambda_{A_2, q_2}g_2, g_1\rangle$$
$$= \langle\mathcal{R}^s_{A_1}u_1, u_2\rangle+ \int_\Omega q_1u_1u_2
- \langle\mathcal{R}^s_{A_2}u_2, u_1\rangle- \int_\Omega q_2u_2u_1$$
$$= \langle\mathcal{R}^s_{A_1}u_1, u_2\rangle- 
\langle\mathcal{R}^s_{A_2}u_1, u_2\rangle
- \int_\Omega (q_2- q_1)u_1u_2$$
\begin{equation}\label{AqId}
=\iint 2(R_{A_2}(x, y)- R_{A_1}(x, y))K(x, y)u_1(y)u_2(x)\,dxdy
- \int_\Omega (q_2- q_1)u_1u_2.
\end{equation}

Now we are ready to prove Theorem 1.1.
\begin{proof}
For $g_j\in C^\infty_c(W_j)$ and $u_j= P_{A_j, q_j}g_j$ ($j= 1,2$),
by (\ref{LIDN}) and (\ref{AqId}) we have
$$\iint G(x, y)u_1(y)u_2(x)\,dxdy
= \int_\Omega (q_2- q_1)u_1u_2$$
where we write
$$G(x, y):= 2(R_{A_2}(x, y)- R_{A_1}(x, y))K(x, y).$$

Note that $\mathrm{supp}\,u_j\subset \Omega\cup W_j$ so the double integral 
on the left hand side is
$$\int_{\Omega\cup {W_2}}\int_{\Omega\cup {W_1}}G(x, y)u_1(y)u_2(x)\,dxdy
= I_1+ I_2+ I_3+ I_4$$
where we write
$$I_1:=\int_\Omega\int_\Omega,\quad I_2:=\int_\Omega\int_{W_1},\quad I_3:=\int_{W_2}\int_\Omega,\quad I_4:=\int_{W_2}\int_{W_1}.$$

Note that $(x, y)\in W_2\times \Omega$ (or $(x, y)\in \Omega\times W_1$) 
implies $(x+ y)/2 \geq r$, $R_{A_1}(x, y)=  R_{A_2}(x, y)=1$, $G(x,y)= 0$ so
$I_2= I_3= 0$. 

Also note that by the assumption on $W^{(1, 2)}$, we can choose $x_0\in W_2, y_0\in W_1$  s.t. $\frac{x_0+ y_0}{2}\notin \mathrm{supp}\,A_j$ so 
$(x, y)\in W_2\times W_1$ implies $(x+y)/2\notin \mathrm{supp}\,A_j$
if we replace $W_2, W_1$ by a small open ball centered at $x_0$ and a small open ball centered at $y_0$ when necessary.
Hence $I_4= 0$ so we have
\begin{equation}\label{OmegaId}
\int_\Omega\int_\Omega G(x, y)u_1(y)u_2(x)\,dxdy
= \int_\Omega (q_2- q_1)u_1u_2.
\end{equation}

Now fix open sets $\Omega_j\subset \Omega$ s.t. 
$\Omega_1\cap \Omega_2= \emptyset$. Also fix
$\phi_j\in C^\infty_c(\Omega_j)$ and $\epsilon > 0$.

By Proposition 2.4, we can choose $g_1\in C^\infty_c(W_1)$ s.t. 
$$||u_1- \phi_1||_{L^2(\Omega)}\leq \epsilon$$
and for this chosen $g_1$, we can choose $g_2\in C^\infty_c(W_2)$ s.t.
$$||u_1||_{L^2(\Omega)}||u_2- \phi_2||_{L^2(\Omega)}\leq \epsilon.$$
Note that $\phi_1(x)\phi_2(x)= 0$ for $x\in \Omega$ so
$$|\int_\Omega (q_2- q_1)u_1u_2|=
|\int_\Omega (q_2- q_1)(u_1-\phi_1)\phi_2+ 
\int_\Omega (q_2- q_1)u_1(u_2- \phi_2)|$$
\begin{equation}\label{se}
\leq ||(q_2- q_1)||_{L^\infty}||\phi_2||_{L^2}||u_1- \phi_1||_{L^2}+
||(q_2- q_1)||_{L^\infty}||u_1||_{L^2}||u_2- \phi_2||_{L^2}\leq C\epsilon.
\end{equation}
Also note that
$$|G(x, y)|\leq 4|\sin{(\frac{x-y}{2}\cdot (A_1- A_2)(\frac{x+y}{2}))}
\sin{(\frac{x-y}{2}\cdot (A_1+ A_2)(\frac{x+y}{2}))}|K(x,y)$$
$$\leq C_A|x- y|^2K(x,y)\leq \frac{C}{|x- y|^{n+ 2s- 2}},$$
which implies
$$\int_\Omega|G(x, y)|dy\leq C_0,\, x\in \Omega,\qquad
\int_\Omega|G(x, y)|dx\leq C_0,\, y\in \Omega.$$
By the generalized Young's Inequality (see Proposition 0.10 on page 9 in \cite{folland1995introduction}),
$$||Tf||_{L^2(\Omega)}\leq C_0||f||_{L^2(\Omega)},\qquad (Tf)(x):= \int_\Omega|G(x, y)f(y)|\,dy$$
so we have
$$|\int_\Omega\int_\Omega G(x, y)u_1(y)u_2(x)\,dxdy
- \int_{\Omega_1}\int_{\Omega_2} G(x, y)\phi_1(y)\phi_2(x)\,dxdy|$$
$$= |\int_\Omega\int_\Omega G(x, y)(u_1(y)- \phi_1(y))\phi_2(x)\,dxdy
+ \int_\Omega\int_\Omega G(x, y)u_1(y)(u_2(x)- \phi_2(x))\,dxdy|$$
$$\leq \int_\Omega\int_\Omega|G(x, y)\phi_2(x)|\,dx|u_1(y)- \phi_1(y)|\,dy
+ \int_\Omega\int_\Omega |G(x, y)u_1(y)|\,dy|u_2(x)- \phi_2(x)|\,dx$$
\begin{equation}\label{de}
\leq C_0||\phi_2||_{L^2}||u_1- \phi_1||_{L^2}
+ C_0||u_1||_{L^2}||u_2- \phi_2||_{L^2}\leq C'\epsilon.
\end{equation}
Combine (\ref{se}), (\ref{de}) with (\ref{OmegaId}). 
Since $\epsilon$ is arbitrary,
$$\int_{\Omega_1}\int_{\Omega_2} G(x, y)\phi_1(y)\phi_2(x)\,dxdy= 0.$$
Note that the set
$$\{\phi_1\otimes \phi_2:\, 
\phi_j\in C^\infty_c(\Omega_j),\,j= 1, 2\}$$
generates a space dense in $C^\infty_c(\Omega_1\times \Omega_2)$
so $G(x, y)= 0$ in $\Omega_1\times \Omega_2$. Since $\Omega_1, \Omega_2$ are arbitrary,
$G(x, y)= 0$ for $x, y \in\Omega$ whenever $x\neq y$ so 
$$R_{A_1}(x, y)= R_{A_2}(x, y),\quad x, y\in\Omega.$$

Now fix $x_0\in \Omega$. Let $A^{(k)}$ denote the $k^{th}$ component of $A$ and let $\{e_k\}^n_{k=1}$ denote
the standard basis of the vector space $\mathbb{R}^n$. Consider 
$x= x_0+ \epsilon e_k$
and $y= x_0- \epsilon e_k$ for small $\epsilon> 0$. Since $|2\epsilon A^{(k)}_j(x_0)|< \frac{\pi}{2}$, $R_{A_1}(x, y)= R_{A_2}(x, y)$ implies 
$|A^{(k)}_1(x_0)|= |A^{(k)}_2(x_0)|$.

Suppose there exist $l\neq k$ s.t. $A^{(k)}_1(x_0)= A^{(k)}_2(x_0)\neq 0$ and
$A^{(l)}_1(x_0)= -A^{(l)}_2(x_0)\neq 0$. Consider
$x= x_0+ \epsilon (e_k+ e_l)$ and $y= x_0- \epsilon (e_k+ e_l)$, then
$$(x-y)\cdot A_j(\frac{x+y}{2})= 2\epsilon (A^{(k)}_j(x_0)+ A^{(l)}_j(x_0)),$$
which contradicts with 
$R_{A_1}(x, y)= R_{A_2}(x, y)$. Hence the only possibility is 
$A_1(x_0)= \pm A_2(x_0)$.
Now we have shown $A_1= \pm A_2$ then Theorem 1.1 is an immediate consequence of Theorem 2.5.
\end{proof}

\section{The Proof of Theorem 1.2}
Now we use the first order linearization in $H^s(\mathbb{R}^n)$ to prove Theorem 1.2 based on Theorem 1.1.
\begin{proof}
For $g\in C^\infty_c(W_1)$ and small $\epsilon> 0$, 
$u^{(j)}_{\epsilon, g}= Q_{A_j, a^{(j)}}(\epsilon g)$ solve
\begin{equation*}
\left\{
\begin{aligned}
\mathcal{R}^s_{A_j} u+ a^{(j)}(x, u)&= 0\quad \text{in}\,\,\Omega\\
u&= \epsilon g\quad \text{in}\,\,\Omega_e\\
\end{aligned}
\right.
\end{equation*}
and $u^{(j)}_g= P_{A_j, a^{(j)}_1}g$ solve
\begin{equation*}
\left\{
\begin{aligned}
\mathcal{R}^s_{A_j} u+ a^{(j)}_1(x)u &= 0\quad \text{in}\,\,\Omega\\
u&=  g\quad \text{in}\,\,\Omega_e\\
\end{aligned}
\right.
\end{equation*}
($j= 1, 2$). By Proposition 2.6, we have
$$u^{(j)}_{\epsilon, g}/\epsilon\to u^{(j)}_g\quad \text{in}\,\,H^s(\mathbb{R}^n),$$
which implies
$$\frac{1}{\epsilon}\mathcal{R}^s_{A_j}u^{(j)}_{\epsilon, g}|_{W_2}\to 
\mathcal{R}^s_{A_j}u^{(j)}_g|_{W_2}\quad\text{in}\,\,
H^{-s}(W_2).$$
Note that (\ref{NIDN}) implies
$$\mathcal{R}^s_{A_1}u^{(1)}_{\epsilon, g}|_{W_2}
=\mathcal{R}^s_{A_2}u^{(2)}_{\epsilon, g}|_{W_2}.$$
Let $\epsilon \to 0$, then we have
$$\mathcal{R}^s_{A_1}u^{(1)}_g|_{W_2}= \mathcal{R}^s_{A_2}u^{(2)}_g|_{W_2},$$
i.e.
$$\Lambda_{A_1, a^{(1)}_1}g|_{W_2}= \Lambda_{A_2, a^{(2)}_1}g|_{W_2},\qquad g\in C^\infty_c(W_1).$$
By Theorem 1.1, $A_1= \pm A_2$. Now Theorem 1.2 is an immediate consequence of Theorem 2.7.
\end{proof}

\bibliographystyle{plain}
\bibliography{Reference2}
\end{document}